\documentclass{article}
\usepackage{graphicx,url}
\usepackage[latin1]{inputenc}
\usepackage{epic,eepic,amsmath,amstext,amssymb,subfigure}

\newtheorem{Def}{Definition}[section]
\newtheorem{lemma}[Def]{Lemma}
\newtheorem{theorem}[Def]{Theorem}

\newtheorem{ex}{Example}

\newcommand{\R}{\mathbb{R}} 
\newcommand{\transpose}{{\hspace{-0.01cm}\scriptscriptstyle\mathrm T}} 
\newcommand{\Oc}{\mathcal O}
\def\eps{\varepsilon}
\def\sfrac{\tfrac}
\def\bigo{\mathcal O}

\title{Numerical integrators for \\ motion under a strong constraining force}

\author{Christian Lubich\footnotemark[2] \ and Daniel Weiss\footnotemark[3]}

\begin{document}
\maketitle

\renewcommand{\thefootnote}{\fnsymbol{footnote}} 
\footnotetext[2]{Mathematisches Institut, Univ.~T\"ubingen, Auf der Morgenstelle,  D-72076 T\"ubingen, Germany. E-mail: {\tt lubich@na.uni-tuebingen.de}}
\footnotetext[3]{Institute for Applied and Numerical Mathematics, Karlsruhe Institute of Technology, Kaiserstr. 89-93, D-76133 Karlsruhe, Germany. E-mail: {\tt daniel.weiss@kit.edu}}

\textbf{Abstract.}
This paper deals with the numerical integration of Hamiltonian systems in which a stiff anharmonic potential causes highly oscillatory solution behavior with solution-dependent frequencies. The impulse method, which uses micro- and macro-steps for the integration of fast and slow parts, respectively, does not work satisfactorily on such problems. Here it is shown that variants of the impulse method with suitable projection preserve the actions as adiabatic invariants and yield accurate approximations, with macro-stepsizes that are not restricted by the stiffness parameter.\\[0.2cm]



\section{Introduction}			
\label{section:introduction}		

We are interested in the efficient numerical integration of Hamiltonian systems in which a stiff anharmonic potential causes highly oscillatory solution behavior with state-dependent slowly varying high frequencies. 

\subsection{The highly oscillatory Hamiltonian system}
We consider Hamiltonians as studied, in varying degrees of generality and with different analytical techniques, by Rubin \& Ungar 
 \cite{RU},  Takens \cite{Tak},  Bornemann \cite{Bor}, Lorenz \cite{L} and Hairer, Lubich \& Wanner \cite[Section XIV.3]{HLW}:
\begin{align}\label{align:HamiltonianGeneral}
 H(x,y) & = \tfrac{1}{2}\,y^{\transpose}M(x)^{-1}y + U(x) + \frac{1}{\varepsilon^{2}}V(x),
 \qquad 0<\varepsilon \ll 1,
\end{align}
depending on positions $x\in\R^n$ and momenta $y\in \R^n$.
The mass matrix $M(x)$ is  symmetric and positive definite and depends smoothly on $x$. The slow potential  $U$ is  smooth, and the stiff potential $\frac{1}{\varepsilon^{2}}V$ with a smooth function $V:D\subset\R^n\to\R$ attains its minimum value $0$ on a $d$-dimensional manifold
\begin{align}\label{align:min}
\mathcal{V} & = \{x\in D \,:\, V(x)=\min V=0\}.
\end{align}
We assume that the potential $V$ is strongly convex along directions non-tangential to~$\mathcal{V}$. More precisely,
there exists $\alpha>0$ such that for $x\in\mathcal{V}$ the Hessian $\nabla^2V(x)$ satisfies
\begin{align}\label{align:quad}
v^\transpose \nabla^2 V(x)v & \geq \alpha \cdot v^\transpose M(x) v
\end{align}
for all vectors $v$ in the $M(x)$-orthogonal complement of the tangent space $T_x \mathcal{V}$. 

Furthermore, we assume that a constraint function $g:D\to \R^m$, with $m=n-d$, is known such that
\begin{align} \label{align:constrain}
\mathcal{V} & = \{x\in D \,:\, g(x)=0\}
\end{align}
and the derivative matrix $G(x)=g'(x)$ is of full rank on $\mathcal{V}$.  

The corresponding system of Hamiltonian differential equations reads
		\begin{align}\label{align:HS}\begin{aligned}
		\dot{x} & = M(x)^{-1}y \\
		\dot{y} & =  -\nabla_x \Big(\tfrac{1}{2}\, y^{\transpose}M(x)^{-1}y\Big) - 
		\nabla U(x) -\frac{1}{\varepsilon^2}\nabla V(x).
			 \end{aligned}\\\nonumber
		\end{align}

\begin{ex}\hspace{0cm}\\
\begin{minipage}{10cm}
A simple, yet nontrivial model example is the stiff spring double pendulum. The Hamiltonian reads
		\begin{align*}
		 H(x,y) & = \tfrac{1}{2}y^{\transpose}y  + U(x) + \frac{1}{\varepsilon^2}V(x),
		\end{align*}
where $\tfrac{1}{2}y^{\transpose}y$ is the kinetic energy, $U(x)=x_{12} + x_{22}$, and 
		\begin{align*}
		\frac{1}{\varepsilon^2}V(x) & =  \tfrac{1}{2}\left[\frac{\alpha_{1}}{\varepsilon}\right]^2\left(\|x_{1}\| - l_1\right)^{2}  + \tfrac{1}{2}\left[\frac{\alpha_{2}}{\varepsilon}\right]^2\left(\|x_{1}-x_{2}\| - l_2\right)^{2}
		\end{align*}
is the stiff potential depending on the small parameter $\varepsilon$. The parameters $l_i$ denote the lengths of the springs, $\alpha_i/\varepsilon$ are the large spring constants.
\end{minipage}\hspace{1.5cm}\begin{minipage}{2cm}
  \includegraphics[scale=0.5]{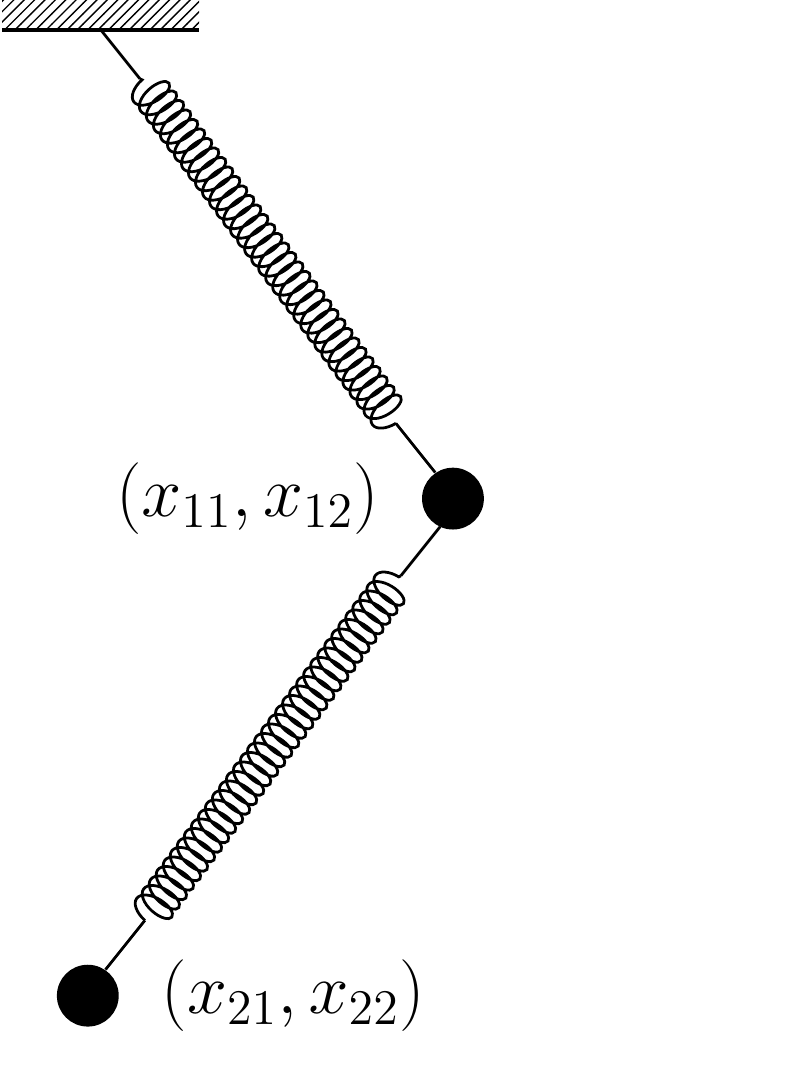}
\end{minipage}\\
\end{ex}

Example 1 helps to fix ideas on a simple toy model. Obviously it extends to chains of stiff springs, which describe the dynamics of chains of atoms in a molecule with almost rigid bonds, cf., e.g., \cite{IRS}.

\subsection{The effective Hamiltonian system}
It has been known since Rubin \& Ungar \cite{RU} that the motion of the system in the limit $\varepsilon\to 0$ differs from the Hamiltonian dynamics constrained to the manifold $\mathcal{V}$ (the rigid double pendulum in the above example) for general initial values $(x_0,y_0)$ that
have an energy bounded independently of~$\varepsilon$,
\begin{align} \label{align:EBound}
H(x_0,y_0) & \leq  \mathit{Const}.
\end{align}
Note that the set of admissible initial values $(x_0,y_0)$ satisfying \eqref{align:EBound} depends on $\varepsilon$.
The effective constrained Hamiltonian has a correction potential $W$,
		\begin{align}\label{align:CH}\begin{aligned}
		 H_{\rm eff}(X,Y) & = \tfrac{1}{2}Y^{\transpose}M(X)^{-1}Y  + U(X) + W(I,X),\\
		 	   0 & = g(X). 
	 	\end{aligned}
		\end{align}
The correction potential $W(I,X)=\sum_{k=1}^m I_k \omega_k(X)$ depends on the $m$ frequencies $\omega_k(X)$, i.e., the square roots of the nonzero generalized eigenvalues of  the pencil $\lambda M(X) -\nabla^2 V(X)$, and on $m$ parameters $I=(I_1,\dots,I_m)$, known as the {\it actions}, which are determined by the initial values $(x_0,y_0)$ of \eqref{align:HS}. The actions vanish for {\it consistent} initial data that satisfy $g(x_0)=0$ and $G(x_0)M(x_0)^{-1}y_0=0$.

As is outlined in \cite{BW} and will be recapitulated in Section \ref{section:transform}, the effective Hamiltonian can be found by transforming the system to separate slow and fast variables as in \cite{L} and \cite[Section XIV.3]{HLW}, and transforming the obtained slow system ($\varepsilon = 0$) back via the effective dynamics of the fast variables.

The effective  constrained Hamiltonian system is then given by
		\begin{align}\label{align:EHS}\begin{aligned}
		\dot{X} & = M(X)^{-1}Y\\
		\dot{Y} & =  -\nabla_X \Big(\tfrac{1}{2}Y^{\transpose}M(X)^{-1}Y\Big) -\nabla U(X) -\nabla_X W(I,X) - G(X)^{\transpose}\Lambda\\
		 	 0 & = g(X),
			 \end{aligned}
		\end{align}
with Lagrange multipliers $\Lambda(t)\in \mathbb{R}^m$. This differs from the usual constrained equations of motion through the correction force $F(I,X)=-\nabla_X W(I,X)$.

To initial values $(x_0,y_0)$  of system \eqref{align:HS} with bounded energy (\ref{align:EBound}),
 we associate consistent initial values $(X_0,Y_0)$  for the effective system  \eqref{align:EHS}. These are chosen by projecting $M$-orthogonally onto the manifold of consistent values:
\begin{align}\label{align:ACIV} \begin{aligned}
X_0 & = x_0+M(x_0)^{-1}G(x_0)^\transpose\lambda,\qquad\;
0  = g(X_0),\\
Y_0 & = y_0+G(X_0)^\transpose\mu,\qquad\qquad\qquad
0  = G(X_0)M(X_0)^{-1}Y_0.
\end{aligned}
\end{align}
With the projection 
\begin{equation}\label{proj}
\mathcal{P}(x)=I-\mathcal{Q}(x), \qquad \mathcal{Q}(x)=[G^\transpose(GM^{-1}G^\transpose)^{-1}GM^{-1}](x),
\end{equation}
 the second equation can be rewritten as $Y_0 = \mathcal{P}(X_0)y_0$, and along the solution of \eqref{align:EHS} we note $Y(t)= \mathcal{P}(X(t))Y(t)$.

 The effective Hamiltonian system describes the limit dynamics on the constraint manifold as long as the solution-dependent frequencies $\omega_k(X(t))$  remain separated and are non-resonant: for some $\delta>0$,
\begin{align}
\label{align:separated}
   \left|\omega_j (X(t))- \omega_k(X(t))\right| &\geq \delta \quad\ \hbox{for } j\ne k
\\
\label{align:nonresonant}
 \left|\omega_j (X(t))\pm \omega_k(X(t)) \pm \omega_l(X(t))\right|& \geq \delta
 \quad\ \hbox{for all $j,\,k,\,l$. }
\end{align}
If these conditions are satisfied for $t\le \overline t$, then we have for the corresponding solutions  of \eqref{align:HS} and  \eqref{align:EHS} over this time interval
\begin{align}  \label{err-eff}
\begin{aligned}
	X(t)-x(t) & = \Oc(\varepsilon) \\
	Y(t) - \mathcal{P}(x(t))y(t) & = \Oc(\varepsilon),
\end{aligned}
\end{align}
where the constants in the $\Oc$-notation depend on $\delta$ and deteriorate as $\delta\to 0$; see \cite{Tak,HLW,BW}. In the case of a single frequency ($m=1$), where no separation and non-resonance conditions appear, the approximation of the full system by the effective system was already studied by  Rubin \& Ungar \cite{RU}. Note that the above estimates also imply 
$$
\mathcal{P}(X(t))(y(t)-Y(t)) = \Oc(\eps),
$$
which is equivalent to saying that the tangential component of the velocity error is $\Oc(\eps)$. The normal component of the velocity is, however, disregarded in the constrained effective equation.

Conditions (\ref{align:separated}) and (\ref{align:nonresonant}) may appear rather severe at first sight, but in fact conditions of this type are needed for the above approximation result for the effective dynamics. Using the techniques of \cite[Chap.\,XIV]{HLW} it can be shown that the order of this approximation is still $\Oc(\varepsilon^{1/(m+1)})$ if $\omega_j \pm \omega_k \pm \omega_l$ have zeros of multiplicity $m$. However, the separation cannot be omitted. If the distance of two frequencies becomes smaller than $\sqrt{\varepsilon}$, then the slow motion can depend very sensitively on the initial values, and it is no longer approximated by the dynamics of the effective Hamiltonian system; see Takens \cite{Tak}.   The indeterminacy of the slow motion in the case of non-separated frequencies is termed {\it Takens chaos} in \cite{Bor}.

\subsection{Outline of the paper and relation to the literature}
The objective of this paper is to devise and analyze a two-scale integrator for the highly oscillatory Hamiltonian system (\ref{align:HS}), such that for a macro-stepsize $h$ that is not restricted by~$\varepsilon$, the method yields an $\Oc(h^2)+\Oc(\varepsilon)$ error in the positions $x(t)$ and the projected momenta $\mathcal{P}(x(t))y(t)$ over time intervals $t=\Oc(1)$.

This paper is part of the vast literature on the numerical solution of highly oscillatory differential equations; see, e.g., the reviews \cite{CJLL,PJY}. Recent work on the numerical integration of highly oscillatory  mechanical systems includes \cite{AST,CSS,SS08,STE,TOM}.

While much work has been done on systems with constant high frequencies, the numerical analysis of the present case of solution-dependent high frequencies or even just the case of explicitly time-dependent high frequencies is scarce; see \cite[Chapter XIV]{HLW}. An important aspect here is to preserve the {\it adiabatic invariants} (see, e.g., \cite{Hen} for this notion) in the numerical discretization.

In this paper we study two-scale time integrators for \eqref{align:HamiltonianGeneral} which aim at solving the effective system \eqref{align:EHS} over the time scale $t=\Oc(1)$ without, however, explicitly evaluating the correction force $F(I,X)=-\nabla_X W(I,X)$. This additional force  is, in general, directly accessible only via a series of computationally expensive, nonlinear implicit coordinate transformations. Moreover, even in cases where the correction force is computationally accessible, it is of interest to have a numerical method  that is able to monitor the possible breakdown of the validity of the effective equation due to the loss of adiabatic invariance of the actions in cases where frequencies come close or become resonant.

Heterogeneous multiscale methods (HMM) \cite{E,ET,EELRV} have been developed for the very purpose to handle  situations where the underlying effective dynamics is not known. In Brumm \& Weiss \cite{BW} an HMM-approach for highly oscillatory mechanical systems with solution-dependent frequencies is analyzed. This approach shows, however, major drawbacks because of difficulties in initializing the micro-simulation. 

In this article, we follow the alternative idea of the impulse method where the Hamiltonian is split into the slow potential $U$ and the fast part including the kinetic and stiff potential energy. The slow part is integrated in macro-steps, the fast part uses micro-steps. As it turns out, this must be complemented with a suitable projection to lead to a method with satisfactory error behavior.

We proceed as follows: In Section \ref{section:Methods} we formulate the impulse method, a mollified impulse method, and a novel projected impulse method for highly oscillatory mechanical systems with solution-dependent frequencies. We state the main convergence theorem and show results of numerical experiments that highlight different behavior of the various methods. In Section \ref{section:transform} we transform the system, following \cite{L} and \cite[Section XIV.3]{HLW} to variables that are appropriate for the further analysis. Moreover, a further mollified impulse method with a projection mollifier in the transformed variables is introduced, which is computationally impractical but serves as a theoretical reference method for the error analysis.  This method is studied in Section~\ref{section:ErrorAna}. Using the obtained results, the analysis of the mollified and projected impulse methods of Section \ref{section:Methods} is done in Section~\ref{section:MPImpulse}.

\section{Numerical methods and statement of the main result}	                 
\label{section:Methods}		                                                                   
\subsection{Impulse method}
The impulse method was introduced in the context of the numerical treatment of molecular dynamics (Grubm\"uller, Heller, Windemuth \& Schulten \cite{GHWS}, Tuckerman, Berne \& Martyna \cite{TBM}). A mathematical study of this method is given by Garc\'{i}a-Archilla, Sanz-Serna \& Skeel \cite{GSS}. The idea is to split the Hamiltonian 
\begin{align*}
H(x,y) & = H^{\mathrm{fast}}(x,y) + U(x)
\end{align*} 
and to approximate the exact flow $\varphi^{H}_{h}$ by the following symmetric decomposition:
\begin{align*}
	\varphi^{H}_{h} \approx \varphi^{\mathrm{slow}}_{h/2}\circ\varphi^{\mathrm{fast}}_{h}\circ\varphi^{\mathrm{slow}}_{h/2}.
\end{align*}
Since the flow of the slow part can be trivially solved, one step is equivalent to 
\begin{enumerate}
\item kick: $y_n^+  = y_n - h/2\cdot \nabla U(x_{n}) $,
\item oscillate: solve system \eqref{align:HS} with $ U = 0$ and initial values  $(y_n^+,x_n)$ over a time step $h$ to obtain $(y_{n+1}^-,x_{n+1})$,
\item another kick: $ y_{n+1}  = y_{n+1}^- -h/2\cdot\nabla U(x_{n+1}) $.
\end{enumerate}
Step 2. is solved approximately using, e.g., the St\"ormer--Verlet method with micro-stepsizes, or alternatively  using a large-timestep method in suitably transformed variables (an adiabatic integrator) as in \cite{L}.
  
Compared to a direct numerical integration of the full system \eqref{align:HS} with small stepsizes, this method saves many evaluations of the slow force $-\nabla U(x)$, which is often the computationally most expensive part.

 \begin{figure}[h]
\begin{center}
\includegraphics[scale = 0.34]{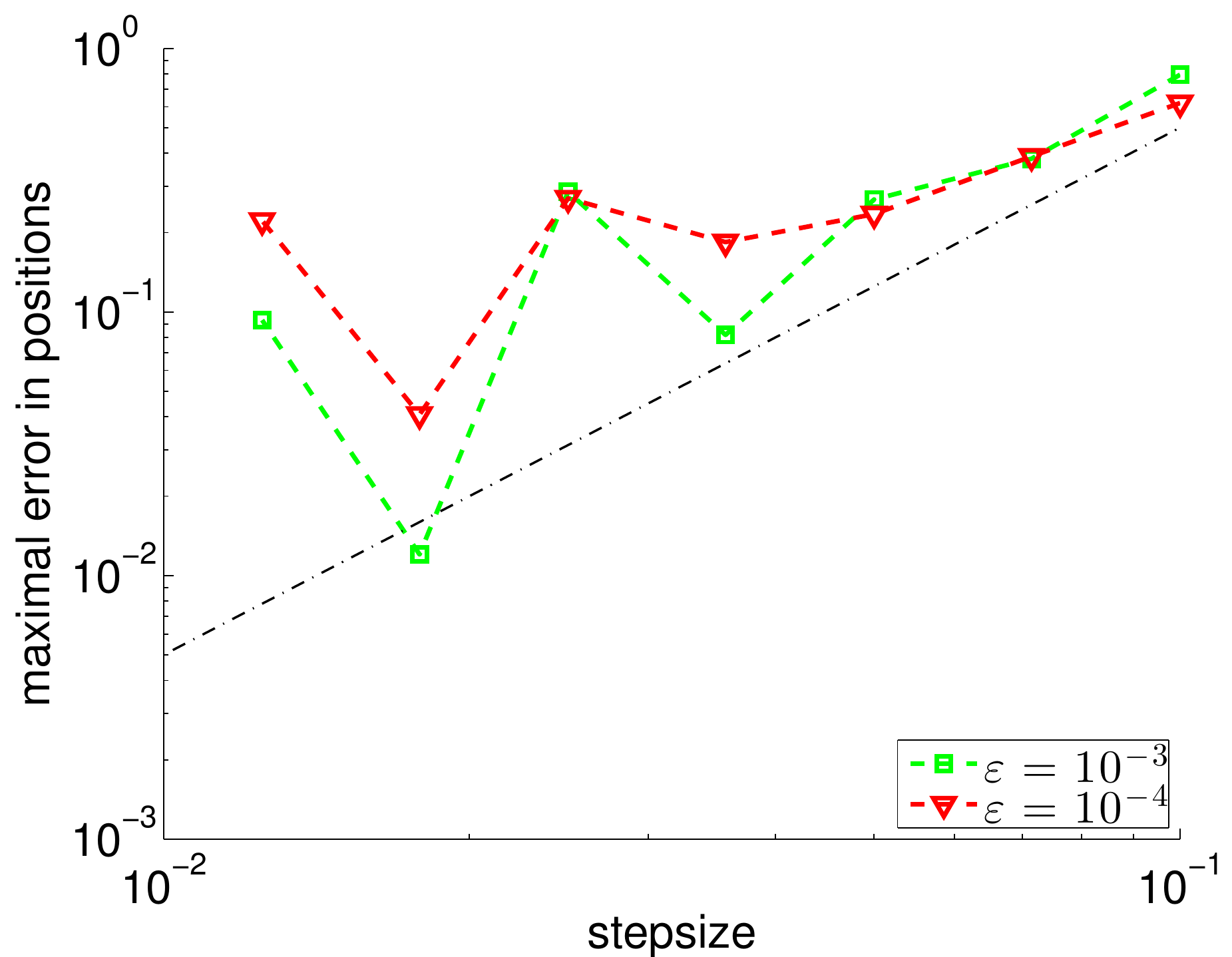}\hspace{1cm}\includegraphics[scale = 0.34]{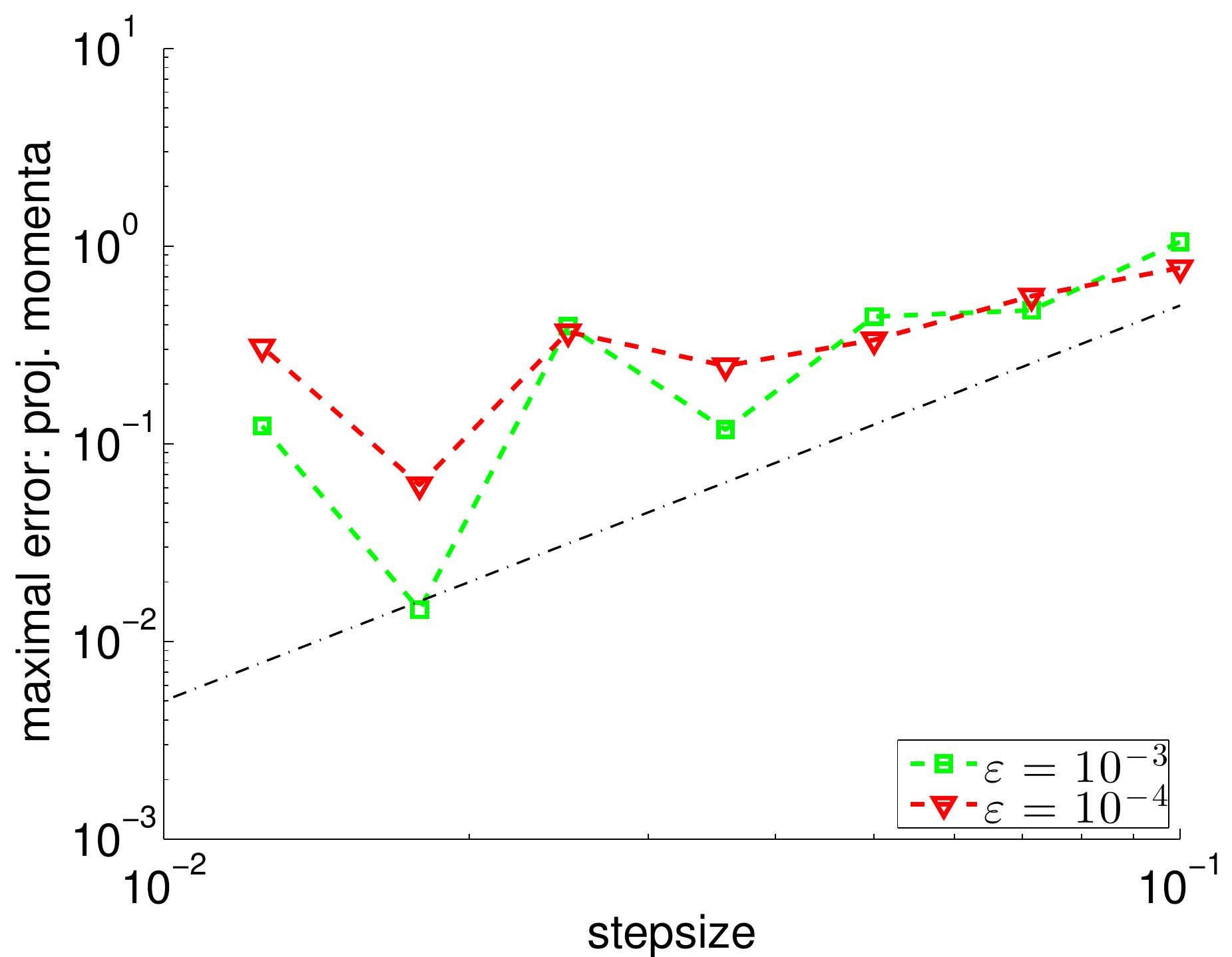}
\end{center}
\caption{Double logarithmic plots: stepsize versus maximal error (maximum norm, maximum over all discrete times) of the impulse method applied to the stiff spring double pendulum with initial value \eqref{align:IV}. Left: Error in positions. Right: Error of the projected momenta.}
\label{figure:IM}
\end{figure}
 
For our numerical experiments we consider the stiff spring double pendulum with initial values
\begin{align} \label{align:IV}
x(0) & = (\sqrt{0.5},-\sqrt{0.5},\sqrt{2},5\varepsilon)^\transpose, & y(0) & = (0,0,0,0)^\transpose
\end{align}
and the parameters $\alpha_1=\alpha_2=1$, over the time interval $0\le t \le 10$. In this situation the frequencies remain well-separated.

We observe unsatisfactory behavior of the impulse method in Figure \ref{figure:IM}. Here and in all following figures the dash-dotted straight line has slope 2, corresponding to the desired $h^2$ error behavior. We used the St\"ormer--Verlet method with very small stepsize ($\eps/1000$ for the impulse method and $\eps/100$ for the following methods) for the micro-integration in order to avoid any significant influence on the overall error. Throughout all computations, as a reference, the effective system \eqref{align:EHS} is approximated in transformed variables (see \eqref{slow-ode}) by the St\"ormer--Verlet method with small stepsize, the results being translated back into cartesian coordinates.

\subsection{Mollified impulse method}
\begin{figure}[h]
\begin{center}
\includegraphics[scale = 0.34]{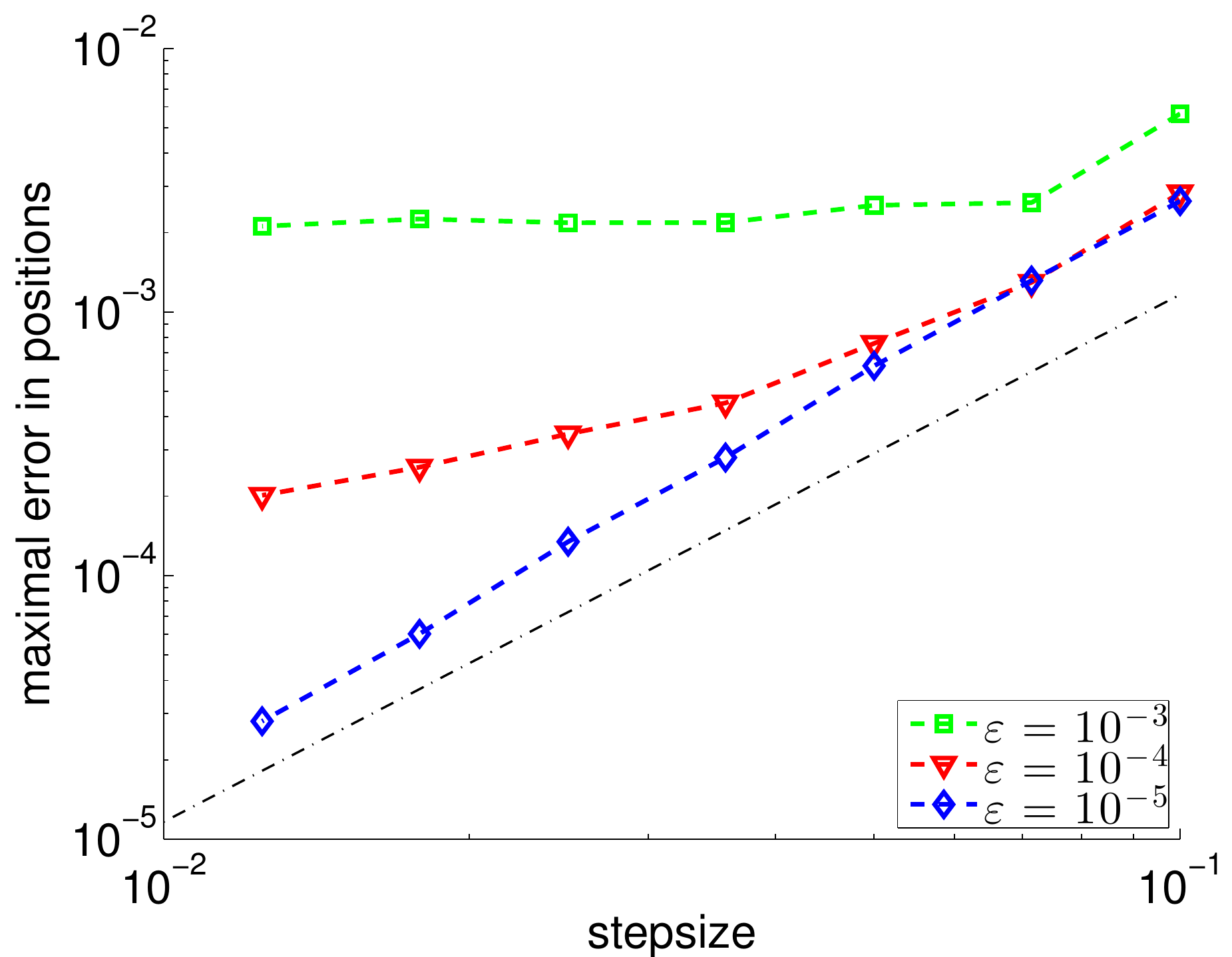}\hspace{1cm}\includegraphics[scale = 0.34]{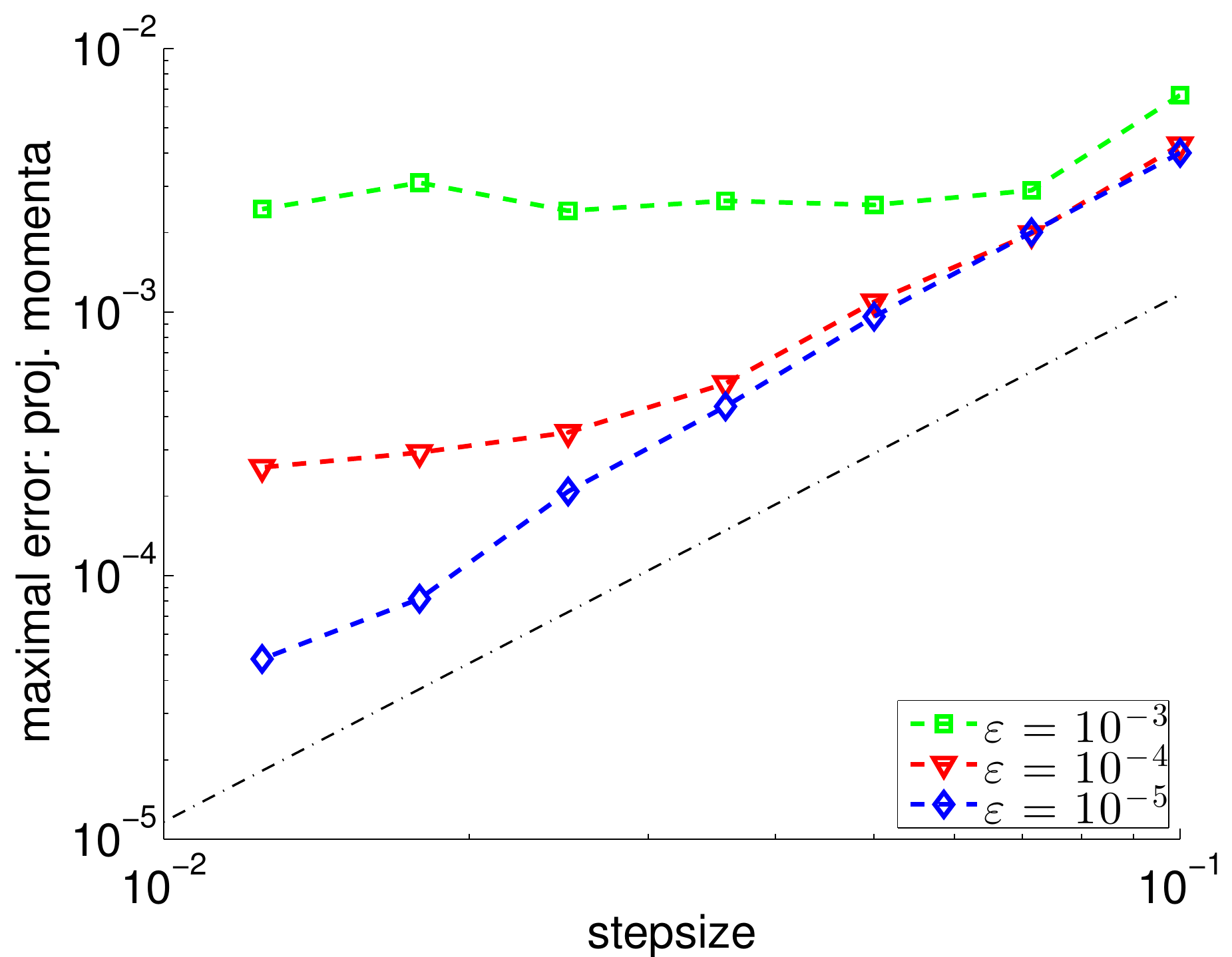}
\end{center}
\caption{Double logarithmic plots: stepsize versus maximal error of the mollified impulse method applied to the stiff spring double pendulum with initial value \eqref{align:IV}. Left: Error in positions. Right: Error of the projected momenta.}
\label{figure:MIM}
\end{figure}

Garc\'{i}a-Archilla, Sanz-Serna \& Skeel \cite{GSS} and Izaguirre, Reich \& Skeel \cite{IRS} improve the impulse method by replacing the slow potential $U(x)$ by a mollified potential $\bar{U}(x)=U(\alpha(x))$, where $\alpha(x)$ is an averaged or suitably projected value of~$x$. The mollified force then reads 
\begin{align*}
 -\nabla \overline{U}(x) & = - \alpha'(x_n)^\transpose\nabla U(\alpha(x_{n})).
\end{align*}
The mollification considered in the present paper is given by the $M(x)$-orthogonal projection onto the configuration manifold $\{X\,:\, g(X)=0\}$, i.e., $\alpha(x)= X$ with
\begin{align}\label{align:ProMol}
\begin{aligned}
X & = x+M(x)^{-1}G(x)^\transpose\lambda,\\
0 & = g(X).
\end{aligned}
\end{align}
Using the same initial value \eqref{align:IV} as in the case of the impulse method, we observe better convergence behavior of the positions and the projected momenta, see Figure~\ref{figure:MIM}. 

As it turns out in the analysis, the unsatisfactory behavior of the impulse method is due to $M(x)$-orthogonal components of the slow forces $\nabla U(x)$. The mollification reduces those $M(x)$-orthogonal components. Indeed, we observe the following.

\begin{lemma} \label{lemma:alpha}
Under the bounded-energy condition $V(x)=\Oc(\eps^2)$, the mollifier $\alpha(x)$ of \eqref{align:ProMol} satisfies
\begin{align*}
\alpha(x) & = x + \Oc(\varepsilon),\\
\alpha'(x)^\transpose &=\mathcal{P}(x)+\Oc(\varepsilon),
\end{align*}
where the projection $\mathcal{P}(x)$ is defined in \eqref{proj}. 
\end{lemma}


\textit{Proof.} The condition $V(x)=\Oc(\eps^2)$ is equivalent to $g(x)=\Oc(\eps)$. Noting that $(GM^{-1}G^\transpose)(x)$ is invertible in view of the full rank of $G$, the implicit function theorem then yields $\lambda=\Oc(\eps)$ such that $g(x+M(x)^{-1}G(x)^\transpose\lambda)=0$,
and hence $\alpha(x)=x+\Oc(\eps)$. Differentiating both equations in \eqref{align:ProMol} yields
\begin{align*}
\alpha'(x) &= I + M^{-1}(x)G(x)^\transpose \lambda'(x) + \Oc(\eps) \\
0 &= G(\alpha(x))\alpha'(x).
\end{align*}
Inserting the first into the second equation permits us to compute 
$$
\lambda'(x)= - (GM^{-1}G^\transpose)^{-1}G(x) + \Oc(\eps).
$$
Reinserting this expression into the first equation yields the stated result on recalling the definition of $\mathcal{P}(x)$.
 
\subsection{Projected impulse method}

The preceding lemma motivates us to simplify the method by projecting  the slow force:
\begin{enumerate}
\item kick: $y_n^+  = y_n - h/2\cdot\mathcal{P}(x_{n})\nabla U(x_{n}) $,
\item oscillate: solve system \eqref{align:HS} with $ U = 0$ and initial values  $(y_n^+,x_n)$ over a time step $h$ obtaining $(y_{n+1}^-,x_{n+1})$,
\item kick: $ y_{n+1}  = y_{n+1}^- -h/2\cdot\mathcal{P}(x_{n+1})\nabla U(x_{n+1}) $.
\end{enumerate}

\begin{figure}[h]
\begin{center}
\includegraphics[scale = 0.34]{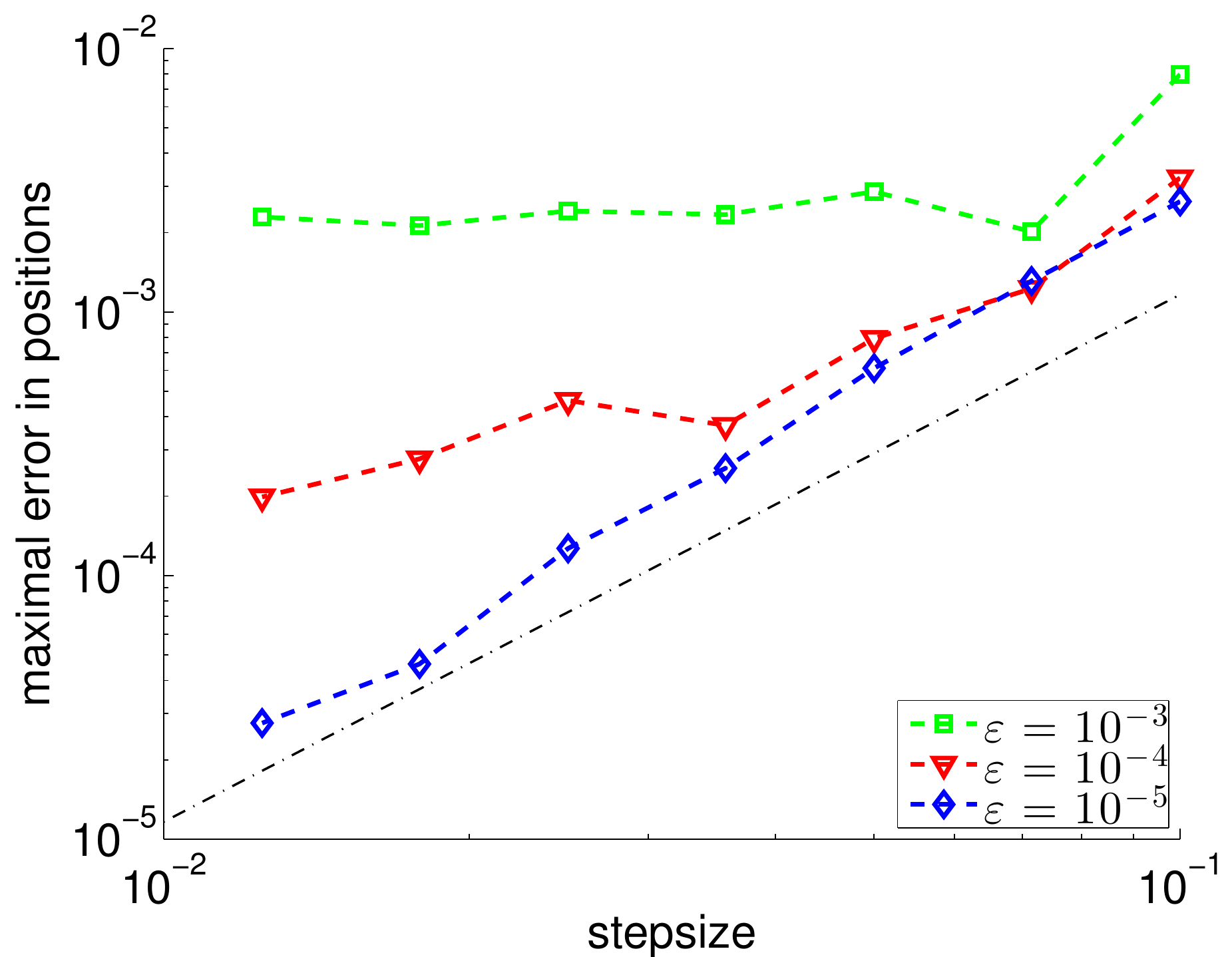}\hspace{1cm}\includegraphics[scale = 0.34]{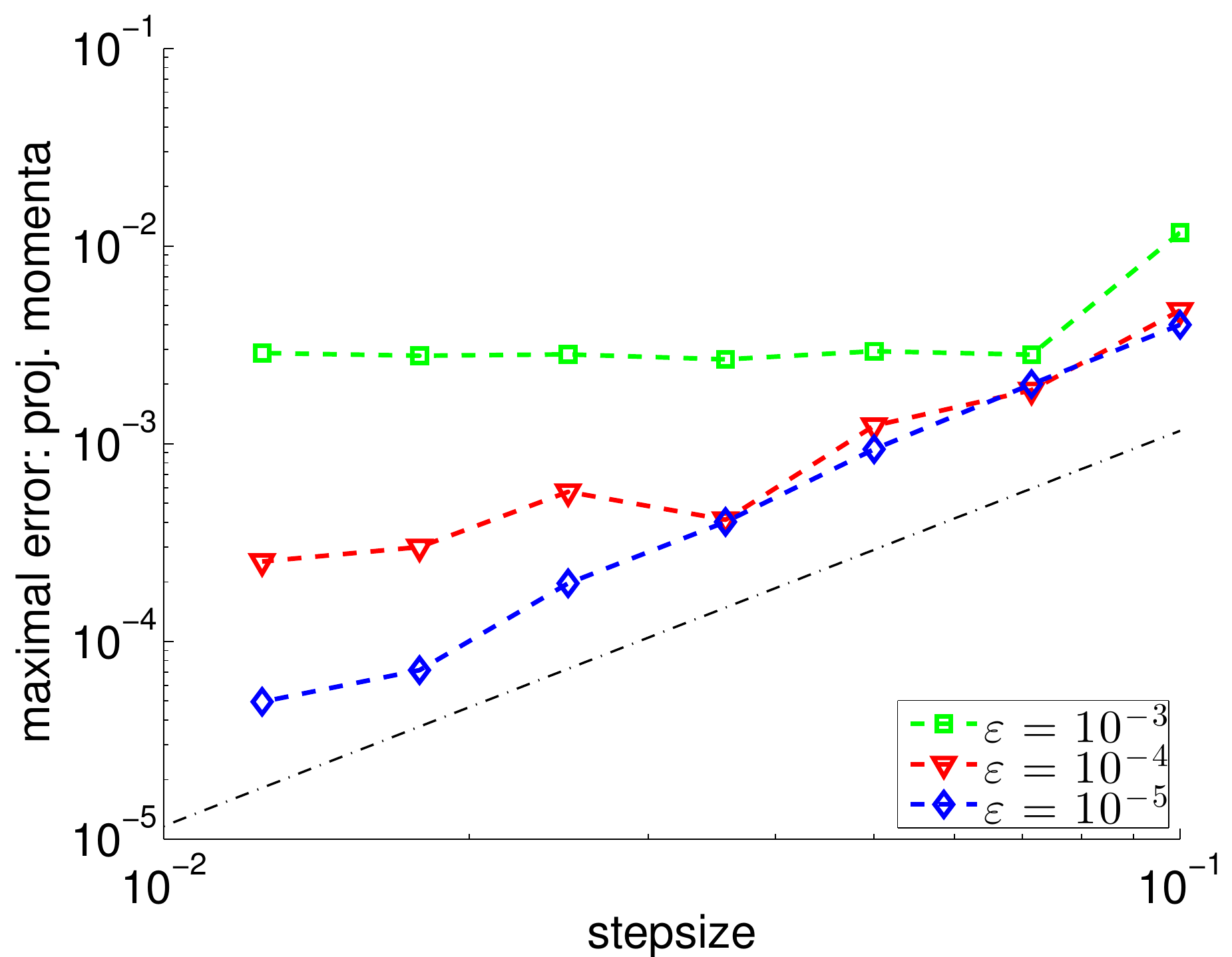}
\end{center}
\caption{Double logarithmic plots: stepsize versus maximal error of the projected impulse method applied to the stiff spring double pendulum with initial value \eqref{align:IV}. Left: Error in positions. Right: Error of the projected momenta.}
\label{figure:PIM}
\end{figure}

Using this new simplified scheme, we observe convergence behavior as in the case of the mollified impulse method, see Figure \ref{figure:PIM}. 

We have sacrificed the symplecticity of the method which is not of main interest here, 
but have nevertheless maintained the time-reversal symmetry.

\subsection{Main Result}
The idea of a projection as a mollification is proposed in \cite{IRS}. There, the use of this idea is shown experimentally but no analysis is given. On the other hand, in \cite{L,HLW} the adiabatic nature of the systems of interest is revealed by applying a series of canonical transformations. Combining the different ideas and techniques,  we are now able to formulate and prove the result about the global error of the mollified impulse method with the projection mollifier. Additionally, we prove the same result for the computationally simpler projected impulse method. The proof is based on a further, different mollification introduced in Section \ref{section:transform}.    


\begin{theorem} \label{theorem:main} Let the initial values satisfy the energy bound \eqref{align:EBound} and assume that the frequencies remain separated and non-resonant (see conditions \eqref{align:separated}-\eqref{align:nonresonant}) along the solution of \eqref{align:HS} for $0\le t \le \overline t$. Then, the errors of the mollified impulse method and of the projected impulse method after $n$ steps with stepsize $h$ satisfy
\begin{align} \label{err-meth-eff}
\begin{aligned}
x_n-X(t_n)& = \Oc(h^2) +\Oc(\varepsilon) \\
\mathcal{P}(x_n)y_n-Y(t_n)& = \Oc(h^2) +\Oc(\varepsilon) 
\end{aligned}
\end{align}
where $(X(t),Y(t))$ is the solution of the effective Hamiltonian system \eqref{align:EHS} with initial values defined by \eqref{align:ACIV}. The constants symbolized by $\Oc$ do not depend on $\varepsilon$, $h$ and $n$ with $nh\leq \overline t$.
\end{theorem}

Combined with \eqref{err-eff} this also yields the error bounds with respect to the solution $(x(t),y(t))$ of the highly oscillatory problem
\begin{align} \label{err-meth-exact}
\begin{aligned}
x_n-x(t_n)& = \Oc(h^2) +\Oc(\varepsilon) \\
\mathcal{P}(x_n)(y_n-y(t_n))& = \Oc(h^2) +\Oc(\varepsilon) .
\end{aligned}
\end{align}
We note, however, that $y_n-Y(t_n) = \Oc(1)$ and $y_n-y(t_n) = \Oc(1)$. Moreover, the method does not converge to the solution $(x(t),y(t))$ of the highly oscillatory system for a fixed $\varepsilon$ as $h\to 0$. This causes no problems since the interest of the method lies in the use of large step sizes $h>\varepsilon$. Note that in Theorem~\ref{theorem:main} there is no restriction of the step size $h$ by the small parameter $\varepsilon$.

Theorem \ref{theorem:main} explains the error behavior observed in Figures 2 and 3.

\section{Transformed variables and another mollified impulse \\ method}	
\label{section:transform}		                                                                   
Under conditions \eqref{align:min}-\eqref{align:constrain},  \cite[Section XIV.3]{HLW} and \cite{L}  show that there exists a canonical change of coordinates $(x,y)=\psi(q,p)$ of the separated form $x=\chi(q)$, $y=\chi'(q)^{-\transpose} p$,
which transforms the Hamiltonian \eqref{align:HamiltonianGeneral} into the form
\begin{eqnarray}
\nonumber
H(q,p) &=& \sfrac12 \, p_0^\transpose  M_0(q_0)^{-1} p_0 + 
\frac1{2\eps}\, p_1^\transpose \Omega(q_0) p_1 +
\frac1{2\eps}\, q_1^\transpose \Omega(q_0) q_1  
\\
\label{eq:osc2-q-ham-rescale}
&&
+\  \sfrac12 \, \begin{pmatrix}
                 p_0 \\ \eps^{-1/2} p_1
                \end{pmatrix}^\transpose
  R(q_0,\eps^{1/2}q_1) \begin{pmatrix}
                 p_0 \\ \eps^{-1/2} p_1
                \end{pmatrix} + \check U(q_0,\eps^{1/2} q_1),
\end{eqnarray}
where $q=(q_0,q_1)\in\mathbb{R}^{d}\times\mathbb{R}^{m}$ and
$p=(p_0,p_1) \in\mathbb{R}^{d}\times\mathbb{R}^{m}$ and the appearing functions have all their partial derivatives bounded independently of $\eps$ and are as follows:
\begin{itemize}
 \item $M_0(q_0)$ is a symmetric positive definite $d\times d$ matrix;
 \item $\Omega(q_0)$ is a diagonal $m\times m$ matrix with positive 
  entries, the {\it frequencies} $\omega_k(q_0)$;
 \item $R(q_0,\eps^{1/2}q_1)$ is a symmetric $n\times n$ matrix with $R(q_0,0)=0$;
 \item $\check U(q_0,\eps^{1/2} q_1) = U(x)$ for $x=\chi(q)$.
\end{itemize}
The assumption (\ref{align:EBound}) of bounded energy  now becomes
\begin{equation}
 \label{qp-bound}
 q=\Oc(\eps^{1/2}), \quad\ p = \Oc(\eps^{1/2}).
\end{equation}
We define the {\it actions}
\begin{equation}
 \label{actions}
 I_k = \frac1{2\eps} \bigl( q_{1,k}^2 + p_{1,k}^2 \bigr), \qquad k=1,\dots, m.
\end{equation}
Under the separation and non-resonance conditions \eqref{align:separated}--\eqref{align:nonresonant}, the actions are {\it adiabatic invariants}: they remain nearly constant along solutions of the Hamiltonian system with bounded energy,
\begin{equation}
 \label{action-inv}
 I_k(t) = I_k(0) + \Oc(\eps),
\end{equation}
see \cite{HLW}, p. 562.

\begin{figure}[h]
\begin{center}
\includegraphics[scale = 0.34]{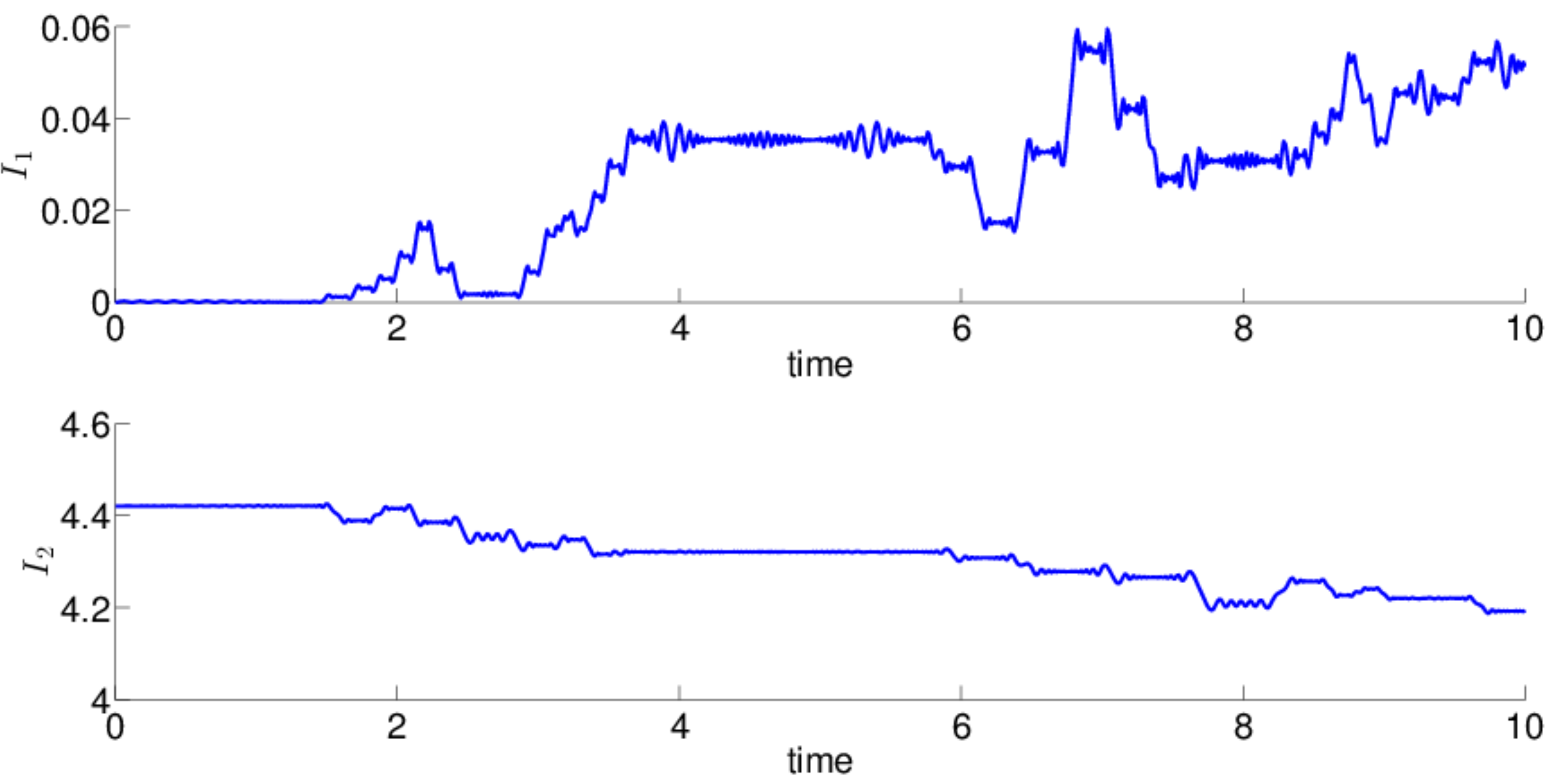}
\end{center}
\caption{Non-preservation of the actions for the impulse method.}
\label{figure:IM_Actions}
\end{figure}

As will become clear from our analysis, this is a key property that should be transfered to the numerical method. However, if we express the impulse method in the transformed variables, then the kick step becomes
$$
\begin{pmatrix}
 p_{n,0}^+ \\ p_{n,1}^+
\end{pmatrix}
= 
\begin{pmatrix}
 p_{n,0} \\ p_{n,1}
\end{pmatrix}
-\frac h2
\begin{pmatrix}
 \nabla_0 \check U (q_{n,0}, \eps^{1/2} q_{n,1}) \\
 \eps^{1/2} \nabla_1 \check U (q_{n,0}, \eps^{1/2} q_{n,1})
\end{pmatrix},
$$
and we see that the actions are not approximately preserved. This is illustrated in Figure~\ref{figure:IM_Actions}. The non-preservation of the actions is at the base of the disappointing numerical behavior observed in Figure~\ref{figure:IM}.
\begin{figure}[h]
\begin{center}
\includegraphics[scale = 0.34]{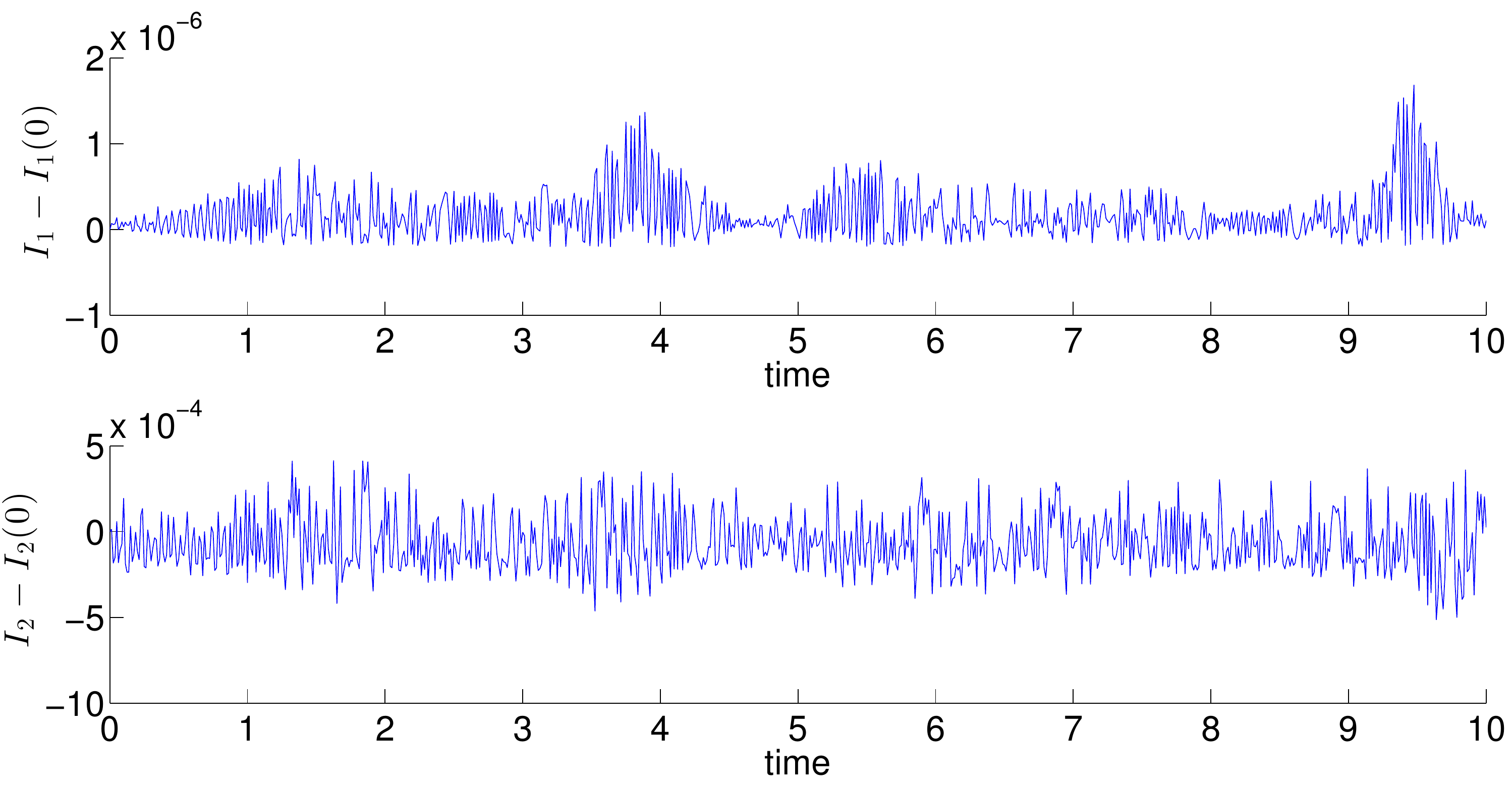}
\end{center}
\caption{Near-preservation of the actions for the projected impulse method.}
\label{figure:PIM_Actions}
\end{figure}

On the other hand, if we use a mollified impulse method for which the modified potential is chosen, in the transformed variables, as
\begin{equation}
\label{moll-transformed}
 \overline U(q) = \check U(q_0,0),
\end{equation}
then $p_{n,1}^+ = p_{n,1}$, and hence the actions are exactly preserved in the kick step. While this method is not practical in that it would require performing the coordinate transformation from $(x,y)$ to $(q,p)$, it gives much theoretical insight into the error propagation behavior. We will therefore study its error in the next section. Subsequently we will interpret the mollified and projected impulse methods of Section 2, which work in the original variables, as perturbations of this theoretically interesting method.

\begin{figure}[h]
\begin{center}
\includegraphics[scale = 0.34]{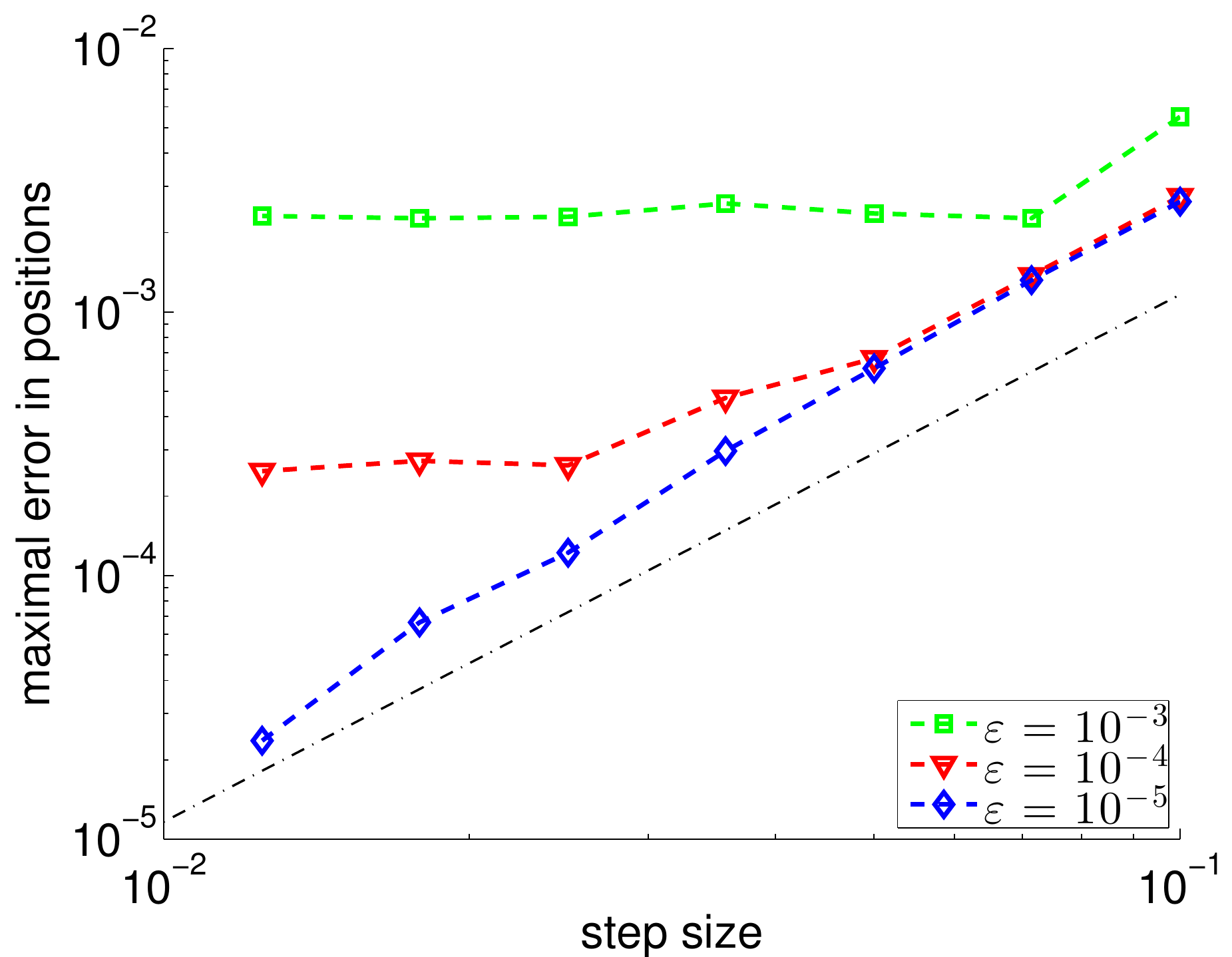}\hspace{1cm}\includegraphics[scale = 0.34]{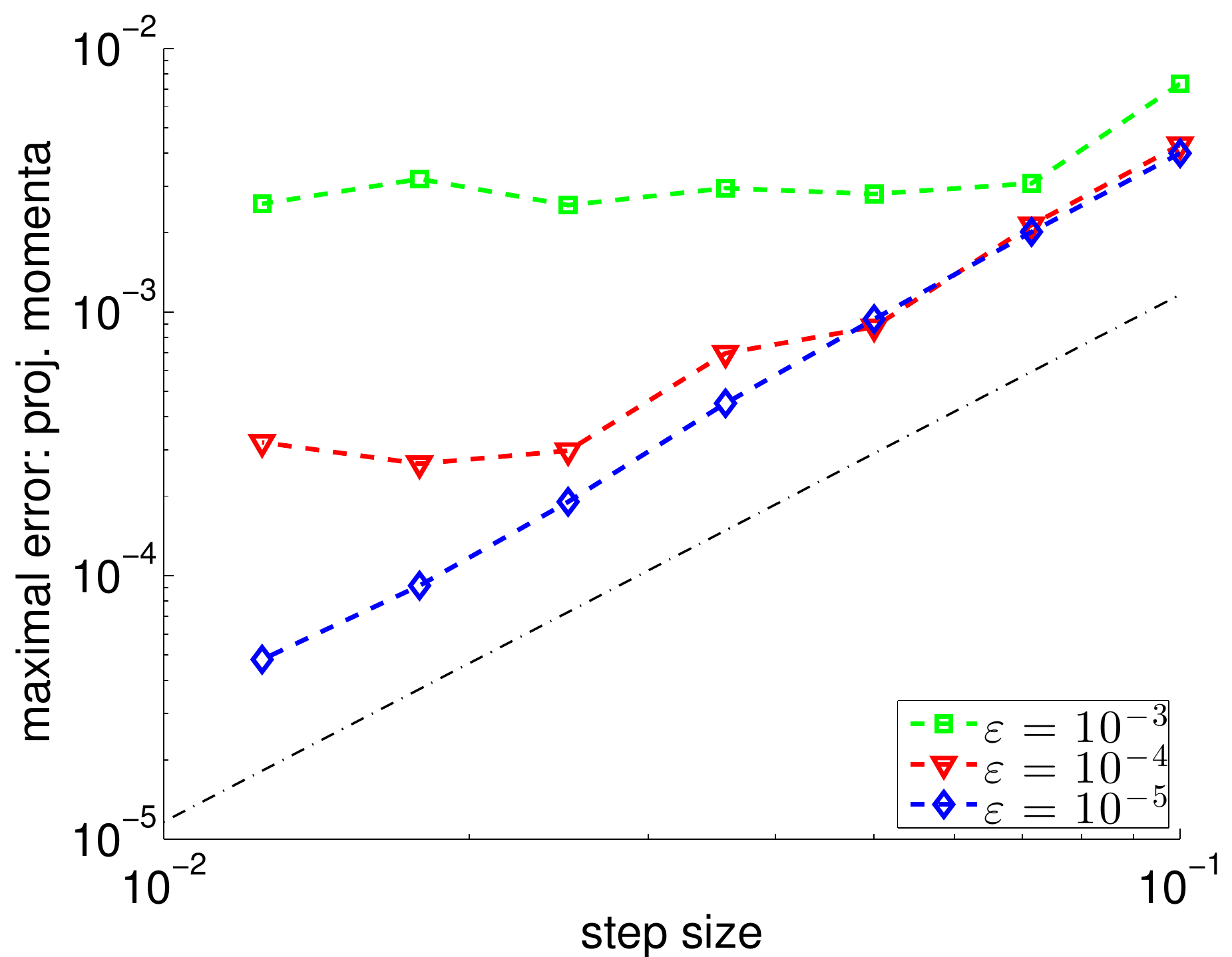}
\end{center}
\caption{Double logarithmic plots: stepsize versus maximal error of the mollified impulse method with \eqref{moll-transformed} applied to the stiff spring double pendulum with initial value \eqref{align:IV}. Left: Error in positions. Right: Error of the projected momenta.}
\label{figure:MIM_TA}
\end{figure}

As a numerical illustration,   in Figure~\ref{figure:MIM_TA} we use the transformed-variable mollified impulse method with the initial value of Section \ref{section:Methods} for the stiff spring double pendulum. We observe similar results as in Figures~\ref{figure:MIM} and \ref{figure:PIM}.  

\section{Error analysis of the transformed-variable method}                              
\label{section:ErrorAna}		                                                                           

We will show almost-conservation of the actions along the numerical solution, $I_n = (I_{n,1}, \dots, I_{n,m})$ with $I_{n,k}= \frac1{2\eps}(q_{n,1,k}^2+p_{n,1,k}^2)$.

\begin{theorem}\label{thm:actions-1}
Assume the energy bound \eqref{align:EBound}. Furthermore, assume that the frequencies of the  transformed-variable mollified impulse method with modified potential \eqref{moll-transformed}
remain separated and non-resonant (see conditions \eqref{align:separated}-\eqref{align:nonresonant})  for $0\le t \le \overline t$. Then, this method approximately preserves the actions:
\begin{align*}
I_n & = I_0 +\Oc(\varepsilon) \quad\ \hbox{ for } nh\leq \overline t.
\end{align*}
The constant symbolized by $\Oc$ is independent of $n$, $h$, and $\eps$.
\end{theorem}

To prove this result, we first need to look in more detail into the differential equations in the transformed variables. As is shown in  \cite{HLW}, p.\,560, the Hamiltonian equations of motion take the form
\begin{eqnarray}
\nonumber
\dot p_0 &=& -\,\nabla_{q_0}
\Bigl(\sfrac12 \, p_0^\transpose  M_0(q_0)^{-1} p_0 + U(q_0,0) \Bigr)
\\
\nonumber
&&-\,\nabla_{q_0}
\Bigl(\frac1{2\eps}\, p_1^\transpose \Omega(q_0) p_1 +
\frac1{2\eps}\, q_1^\transpose \Omega(q_0) q_1 \Bigr) + f_0(p,q)
\\[1mm]
\label{pq-ode}
\dot q_0 &=& M_0(q_0)^{-1} p_0 + g_0(p,q)
\\[2mm] 
\nonumber
\begin{pmatrix} \dot p_1 \\ \dot q_1 \end{pmatrix} &=&
\frac1\eps
\begin{pmatrix}0 & -\Omega(q_0) \\ \Omega(q_0) & 0\end{pmatrix}
\begin{pmatrix} p_1 \\  q_1\end{pmatrix} +
\begin{pmatrix} f_1(p,q) \\ g_1(p,q) \end{pmatrix}
\end{eqnarray}
with  functions
$f_0=\bigo(\eps)$, $g_0=\bigo(\eps)$ and
$f_1=\bigo(\eps^{1/2})$, $g_1=\bigo(\eps^{1/2})$.
Moreover we have (omitting the arguments $p_0,q_0$ in 
$a,b,c,L$, which are all $\bigo(1)$)
\begin{eqnarray}
\nonumber
f_1 &=& -\eps^{1/2}c - Lp_1 + \eps^{-1/2} a(p_1,p_1)
- \eps^{1/2}  \nabla_1 \check U(q_0,0) + \bigo(\eps^{3/2})
\\
\label{f1g1}
g_1 &=& \qquad \quad\ \ \, L^T q_1 + \eps^{-1/2} b(p_1,q_1)
 + \bigo(\eps^{3/2})
\end{eqnarray}
where $L$ is an $m\times m$ matrix and the functions $a$ and $b$ are bilinear.

We diagonalize
$$
\Gamma^*\begin{pmatrix}0 & -\Omega(q_0) \\ \Omega(q_0) & 0\end{pmatrix} \Gamma = i\begin{pmatrix} \Omega( {q}_0) & 0 \\ 0 & -\Omega( {q}_0)\end{pmatrix}=: i\Lambda(q_0),
$$
with 
$$
\Gamma  = \frac{1}{\sqrt{2}}\begin{pmatrix} I & I \\ -iI & iI \end{pmatrix},
$$
and introduce the diagonal phase matrix $\Phi(t)$  by
\begin{align*}
{\Phi}(t) & = \int_0^t \Lambda( {q}_0(s))\, ds.
\end{align*} 
Following \cite{HLW}, p.\,561, we transform the oscillatory part of the solution to {\it adiabatic variables}
\begin{align}\label{align:adiabaticTrans}
\eta(t) & = \varepsilon^{-1/2}\exp \left(-\frac{i}{\varepsilon}\Phi(t)\right)\Gamma^* \begin{pmatrix}  {p}_1(t)\\ {q}_1(t)\end{pmatrix}.
\end{align}  
We further introduce the $m\times 2m$ matrices $P_1(t)$ and $Q_1(t)$ as
\begin{align*}
 \begin{pmatrix} P_1\\Q_1\end{pmatrix}& =\Gamma\exp\left(\frac{i}{\varepsilon}\Phi\right),
\end{align*}
so that $p_1=\eps^{1/2}P_1\eta$ and $q_1=\eps^{1/2}Q_1\eta$. In adiabatic variables, the differential equation for the oscillatory part becomes
\begin{eqnarray}
\nonumber
\dot \eta &=& \exp\Bigl( -\frac i\eps \Phi\Bigr) W(p_0,q_0) 
\exp\Bigl( \frac i\eps \Phi\Bigr)\eta 
\\
\label{eq:osc2-q-eta-ode-2}
&&
+\
\exp\Bigl( -\frac i\eps \Phi\Bigr) \,\Gamma ^* 
\begin{pmatrix}
 a(P_1\eta,P_1\eta;p_0,q_0) \\
 b(P_1\eta,Q_1\eta;p_0,q_0) \end{pmatrix}
\\
\nonumber
&&
-\  P_1^*\Bigl( c(p_0,q_0) +  \nabla_1\check U(q_0,0) \Bigr)
+ \bigo(\eps)\end{eqnarray}
with
\begin{align*}
W & = -\tfrac{1}{2}\begin{pmatrix}L-L^\transpose & L+L^\transpose \\ L+L^\transpose & L-L^\transpose \end{pmatrix}.
\end{align*}
The functions $L,a,b,c$ are those appearing in the remainder terms $f_1$ and $g_1$ in \eqref{f1g1}.

\textit{Proof.} (of Theorem~\ref{thm:actions-1})  We rewrite the mollified impulse method in adiabatic variables and note that the kick steps do not change the adiabatic variables: in the $j$th time step, $\eta_j^+=\eta_j$ and $\eta_{j+1}=\eta_{j+1}^-$. We thus obtain
\begin{align*}
\eta_{j+1} & =  \eta_j  + \int_{t_j}^{t_{j+1}} \dot{\eta}^j(s)ds,\\
\end{align*}
where $\eta^j(t)$ solves 
\begin{align}\label{align:Dynamicseta}\begin{aligned}
 \dot{\eta} & = \exp\left(-\frac{i}{\varepsilon}\Phi^j\right) W( {p}_0^j, {q}_0^j) \exp\left(\frac{i}{\varepsilon}\Phi^j\right)\eta\\
 		 &\quad + \exp\left(-\frac{i}{\varepsilon}\Phi^j\right)\Gamma^*\begin{pmatrix} a(P_1^j\eta,P_1^j\eta; {p}_0^j, {q}_0^j) \\ b(P_1^j\eta,Q_1^j\eta; {p}_0^j, {q}_0^j) \end{pmatrix} + {(P_1^j)}^* c( {p}_0^j, {q}_0^j) +\Oc(\varepsilon)
\end{aligned}		 
\end{align}
with initial value  $\eta_j$ on the interval $[t_j,t_{j+1}]$. All terms with superscript $j$ are defined with respect to the solution of the oscillation step of the mollified impulse method on $[t_{j},t_{j+1}]$. 

In the remaining part of the proof we show $\sum_{j=0}^{n-1} \int_{t_j}^{t_{j+1}} \dot{\eta}^j(s)ds=\Oc(\varepsilon)$. The techniques are more or less the same as for the exact solution presented in \cite{HLW}. Therefore, we just consider the first term of the righthand side in \eqref{align:Dynamicseta}. For $l\not=k$ partial integration gives
\begin{align*}
\sum_{j=0}^{n-1}\int_{t_j}^{t_{j+1}} \exp\left(-\frac{i}{\varepsilon}(\Phi_l^j(s)-\Phi_k^j(s))\right)w_{lk}(p_0^j(s),q_0^j(s))\eta_k^j(s) ds & =  \\
& \hspace{-7cm}i\varepsilon\sum_{j=0}^{n-1}  \exp\left(-\frac{i}{\varepsilon}(\Phi_l^j(s)-\Phi_k^j(s))\right)\frac{w_{lk}(p_0^j(s),q_0^j(s))\eta_k^j(s)}{\omega_l^j(q_0(s))-\omega_k^j(q_0(s))}  \Big|_{t_j}^{t_{j+1}}\\
& \hspace{-7cm}-i\varepsilon\sum_{j=0}^{n-1} \int_{t_j}^{t_{j+1}} \exp\left(-\frac{i}{\varepsilon}(\Phi_l^j(s)-\Phi_k^j(s))\right)\frac{d}{ds}\frac{w_{lk}(p_0^j(s),q_0^j(s))\eta_k^j(s)}{\omega_l^j(q_0(s))-\omega_k^j(q_0(s))} ds,
\end{align*}
where the latter term is of size $\Oc(\varepsilon)$ in the case of separated frequencies. Taking into account the $\Oc(h)$-jumps from $p_0^{j}(t_{j+1})$ to $p_0^{j+1}(t_{j+1})$ and noting that $\eta^{j}(t_{j+1})=\eta^{j+1}(t_{j+1})$ and $\Phi^{j}(t_{j+1})=\Phi^{j+1}(t_{j+1})$, we prove the same bound for the first term.
We have thus shown 
\begin{equation}\label{eta-inv}
 \eta_n=\eta_0+ \bigo(\eps),
\end{equation}
and since $I_{n,k}=|\eta_{n,k}|^2$, the result follows.\\[0.5cm]

We are now in the situation to prove an error bound.

\begin{theorem} \label{theorem:help} Assume the energy bound \eqref{align:EBound}. Furthermore, assume that the frequencies remain separated and non-resonant (see conditions \eqref{align:separated}-\eqref{align:nonresonant}) on a fixed time interval $0\le t \le \overline t$. Then, the error of the transformed-variable mollified impulse method of Section \ref{section:transform} after $n$ steps with stepsize $h$ satisfies
\begin{align*}
x_n-X(t_n)& = \Oc(h^2) +\Oc(\varepsilon)\\
\mathcal{P}(x_n)y_n-Y(t_n)& = \Oc(h^2) +\Oc(\varepsilon) .
\end{align*}
The constants symbolized by $\Oc$ do not depend on $\varepsilon$, $h$ and $n$ with $nh\leq \overline t$.
\end{theorem}

We note, however, that the normal components of the momenta are not approximated correctly: we only have $y_n-y(t_n) = \Oc(1)$.

\textit{Proof.} We consider the method in the slow components $p_0,q_0$ as a perturbed variant of the St\"ormer--Verlet scheme applied to the slow system
\begin{align}
\label{slow-ode}
\begin{aligned}
\dot{p}_0 & = -\nabla_{q_0}\left(\tfrac{1}{2}p_0^\transpose M_0(q_0)p_0 + U(q_0,0)\right)-\sum_{k=1}^m I_k(0)\nabla_{q_0}\omega_k(q_0),\\
\dot{q}_0 & =  M_0(q_0)^{-1}p_0.
\end{aligned}
\end{align}
More precisely, if we write the St\"ormer--Verlet scheme for \eqref{slow-ode} in one-step form  as
\begin{align*}
\begin{pmatrix} p_{n+1,0}\\q_{n+1,0} \end{pmatrix} & = \Psi_h(p_{n,0},q_{n,0}),
\end{align*}
then the slow components of the mollified impulse method for \eqref{pq-ode} fulfill 
 \begin{align*}
\begin{pmatrix} p_{n+1,0}\\q_{n+1,0} \end{pmatrix} & = \Psi_h(p_{n,0},q_{n,0}) + d_n
\end{align*}
with a local error $d_n=\Oc(h^3)+\Oc(h\varepsilon)$, because 
$$
\nabla_{q_0}
\Bigl(\frac1{2\eps}\, p_1^\transpose \Omega(q_0) p_1 +
\frac1{2\eps}\, q_1^\transpose \Omega(q_0) q_1 \Bigr) = \sum_{k=1}^m I_k\nabla_{q_0}\omega_k(q_0)
$$
and $I_{n,k}=I_{0,k}+\bigo(\eps)$ by Theorem~\ref{thm:actions-1}.
Application of the discrete Gronwall Lemma gives the desired result for the slow components in the variables $(q,p)$: for $nh\le\overline t$,
$$
q_{n,0}=q_0(t_n)+\bigo({h^2}) +\bigo(\eps), \quad\
p_{n,0}=p_0(t_n)+\bigo({h^2}) +\bigo(\eps).
$$
For the fast variables we have, using (\ref{eta-inv}),
\begin{align*}
\begin{pmatrix} p_1(t_n) \\ q_1(t_n)\end{pmatrix} & = \varepsilon^{1/2}\Gamma\exp\left(\frac{i}{\varepsilon}\Phi(t_n)\right) \eta(t_n) \\
& = \varepsilon^{1/2}\Gamma\exp\left(\frac{i}{\varepsilon}\Phi(t_n)\right) [\eta(0)+\Oc(\varepsilon)]\\
& = \varepsilon^{1/2}\Gamma\exp\left(\frac{i}{\varepsilon}\Phi(t_n)\right) [\eta_n +\Oc(\varepsilon)]\\
& = \Gamma\exp\left(\frac{i}{\varepsilon}[\Phi(t_n)-\Phi^n(t_n)]\right)\Gamma^* \begin{pmatrix}p_{n,1}\\q_{n,1} \end{pmatrix}+\Oc(\varepsilon^{3/2}),
\end{align*}  
which shows 
$$
q_{n,1}=\bigo(\eps^{1/2}),\quad \ p_{n,1}=\bigo(\eps^{1/2}),
$$
but in view of the phase difference $\Phi(t_n)-\Phi^n(t_n)=\bigo(h^2)+\bigo(\eps)$ this does not yield an approximation estimate.
Transforming back to the coordinates $(x,y)$ and considering the rescaling and the Lipschitz-continuity of the transformations then gives the result.

\section{Error analysis of the mollified and projected impulse methods}	
\label{section:MPImpulse}		                                                                   
In order to analyze the mollified  and projected impulse methods of Section \ref{section:Methods}, we have to derive an appropriate expression of the kick-step in the transformed variables $(q,p)$. We show that both methods are $\Oc(\varepsilon)$-perturbations of the transformed-variable mollified impulse method of Section~\ref{section:transform}, which uses the modified potential \eqref{moll-transformed}, that is, in the original variables, the modified potential  $U(\pi(x))$ with
$$
\pi(x) = \chi(q_0,0) \quad\hbox{for } x=\chi(q) \hbox{ with } q=(q_0,q_1).
$$
The following result is essential in relating the various methods.

\begin{lemma} \label{lemma:pi}
For the mollifier $\pi(x)$ we have, under the bounded-energy condition $V(x)=\Oc(\eps^2)$,
\begin{align*}
\pi(x) & = x + \Oc(\varepsilon),\\
\pi'(x)^\transpose &=\mathcal{P}(x)+\Oc(\varepsilon),
\end{align*}
where the projection $\mathcal{P}(x)$ is defined in \eqref{proj}.
\end{lemma}

\textit{Proof.} As the construction in \cite[Chapter XIV.3]{HLW} shows, the transformation $x=\chi(q)$ is composed as $\chi=\xi\circ\phi_\eps$ with an $\eps$-independent transformation $\xi$ and a rescaling $\phi_\eps(q_0,q_1)=(q_0,\eps^{1/2}q_1)$. Since for $x=\chi(q)$ the bounded-energy condition $V(x)=\Oc(\eps^2)$ is equivalent to $q_1=\Oc(\eps^{1/2})$, we obtain
$$
x=\xi(q_0,\eps^{1/2} q_1) = \xi(q_0,0)+ \Oc(\eps) = \pi(x)+\Oc(\eps).
$$
The proof of the result for $\pi'(x)^\transpose$ is then obtained from the identity
$$
\pi'(X)^\transpose = \mathcal{P}(X) \quad\hbox{ for $X$ with $g(X)=0$}
$$
used for $X=\pi(x)$. This identity is obtained from the transformation laws as follows: Under a change of coordinates $x=\chi(q)$, the constraint function changes to $\check g(q)=g(\chi(q))$ and its derivative $\check G=\check g'$ to $\check G(q) = G(x) \chi'(q)$. The inverse mass matrix changes to 
$\check M(q)^{-1} = \chi'(q)^{-1} M(x) ^{-1}\chi'(q)^{-\transpose}$. Consequently, the projection $\check {\mathcal{P}} (q) = I - [\check G^\transpose (\check G \check M^{-1} \check G^\transpose)^{-1} \check G \check M^{-1} ](q)$ transforms as
$$
\check { \mathcal{P} }(q) = \chi'(q)^\transpose \mathcal{P}(x) \chi'(q)^{-\transpose}.
$$
For $\check\pi(q) = \chi^{-1}(\pi(\chi(q)))$, the transposed derivative transforms in the same way for $x=\chi(q)$ with $x=\pi(x)$:
$$
\check\pi'(q)^\transpose = \chi'(q)^\transpose \pi'(x)^\transpose \chi'(q)^{-\transpose}.
$$
Now, in the variables $q$ we have $\check g(q_0,q_1)=q_1$ and a block diagonal mass matrix $\check M(q)$, and on the other hand $\check\pi(q_0,q_1)=(q_0,0)$.  This gives us
$$
\check {\mathcal{P}}(q) = \begin{pmatrix} 1 & 0 \\ 0 & 0\end{pmatrix} = \check \pi'(q)^\transpose,
$$
and hence the result follows. \\[0.5cm]

Using Lemma \ref{lemma:pi} and the corresponding result Lemma \ref{lemma:alpha} for the mollified impulse method of Section~\ref{section:Methods}, we find for the kick step of the mollified - and projected impulse methods expressed in the variables $(q,p)$ of Section~\ref{section:transform}
\begin{align*}
\begin{pmatrix}  {p}_{n,0}^{+} \\  {p}_{n,1}^{+} \end{pmatrix} & = \begin{pmatrix}  {p}_{n,0} \\  {p}_{n,1} \end{pmatrix} -\frac{h}{2}\begin{pmatrix} \nabla_{q_{0}}\check U( {q}_{n,0}, 0) \\ 0 \end{pmatrix} + \begin{pmatrix} \Oc(h\varepsilon) \\ \Oc(h\varepsilon^{3/2}) \end{pmatrix}. 
\end{align*}
Here, the additional factor $\varepsilon^{1/2}$ is due to the rescaling of the fast positions and momenta in the transformation. For the actions, we then obtain the estimate  $ I_n^{+}  = I_{n} +\Oc(h\varepsilon)$ in the kick step, and as in 
Section~\ref{section:ErrorAna} we obtain the near-invariance of the actions along the numerical solution.

\begin{theorem}\label{thm:actions-2} Consider the mollified impulse method or the projected impulse method of Section~\ref{section:Methods}.
Assume the energy bound \eqref{align:EBound}. Furthermore, suppose that the frequencies 
remain separated and non-resonant (see conditions \eqref{align:separated}-\eqref{align:nonresonant}) for $0\le t \le \overline t$. Then, the method approximately preserves the actions:
\begin{align*}
I_n & = I_0 +\Oc(\varepsilon) \quad\ \hbox{ for } nh\leq \overline t.
\end{align*}
The constant symbolized by $\Oc$ is independent of $n$, $h$, and $\eps$.
\end{theorem}

Therefore, Theorem \ref{theorem:main} holds true, which can be proven in exactly the same way as Theorem \ref{theorem:help}.

\section{Conclusion}
\label{section:conclusion}
We devised and analyzed numerical integrators that capture the effective dynamics of stiff mechanical systems where a strong constraining force leads to highly oscillatory solution behavior with state-dependent frequencies. The integrators are projected variants of the impulse method, which splits the Hamiltonian into the slow potential and the fast part and integrates the latter with micro-steps. A key aspect is the preservation of the actions as adiabatic invariants in the numerical method. It is shown that this can be ensured by 
approximately projecting out the normal components of the slow force. Numerical experiments and theoretical error bounds illustrate the favorable properties of the proposed methods.

\bibliography{Literatur}

\end{document}